\documentclass[a4paper,12pt]{article}
\usepackage{tipa,amsmath,amsfonts,amssymb,amsthm,amscd,mathrsfs,latexsym,bm}
\usepackage[all]{xy}
\usepackage{graphicx}
\usepackage{mathptm,pslatex}

 \newcommand{\lam}{\lambda}

 \DeclareMathOperator{\rad}{rad}

\def\qed{\hfill$\Box$\smallskip}

\newtheorem{prop}{Proposition}[section]
\newtheorem{thm}[prop]{Theorem}\newtheorem{cor}[prop]{Corollary}
\newtheorem{lem}[prop]{Lemma}\newtheorem{dfn}[prop]{Definition}

\title{Centers of symmetric cellular algebras
\thanks {keywords: symmetric cellular algebras, center.}}
\date{November 18, 2009}

\author{Yanbo Li\\[5pt]Department of Information and Computing Sciences,
\\Northeastern University at Qinhuangdao; \\Qinhuangdao, 066004, P.R. China  \\[5pt]School of Mathematics Sciences,
Beijing Normal University;\\
Beijing, 100875, P.R. China\\[1.5pt] E-mail: liyanbo707@163.com}
\begin{document}
\maketitle

\begin{abstract}

Let $R$ be an integral domain and $A$ a symmetric cellular
algebra over $R$ with a cellular basis $\{C_{S,T}^\lam \mid \lam\in\Lambda, S,T\in
M(\lam)\}$. We will construct an ideal $L(A)$ of the center of $A$ and prove that $L(A)$ contains the so-called Higman ideal. When $R$ is a field, we prove that the dimension of $L(A)$ is not less than the number of non-isomorphic simple $A$-modules.
\end{abstract}

\section{Introduction}

In 1996, Graham and Lehrer \cite{GL} introduced cellular algebras
in order to provide a systematic
framework for studying the representation theory of a class of algebras. By the theory of cellular algebras, one can
parameterize simple modules for a finite dimensional cellular
algebra by methods in linear algebra. Many classes of algebras from
mathematics and physics are found to be cellular, including Hecke
algebras of finite type, Ariki-Koike algebras, $q$-Schur algebras,
Brauer algebras, Temperley-Lieb algebras, cyclotomic Temperley-Lieb
algebras, partition algebras, Birman-Wenzl algebras
and so on, see \cite{G}, \cite{GL}, \cite{RX}, \cite{Xi1},
\cite{Xi2} for details.

There are many papers on centers of Hecke algebras of finite type,
which are all cellular algebras \cite{G}. In \cite{J},  Jones found bases for centers of Hecke algebras of
type A over $\mathbb{Q}[q,q^{-1}]$, where $q$ is an indeterminant.
This basis is an analog of conjugacy class sum in a group algebra.
In \cite{GR}, Geck and Rouquier found bases for the centers of
generic Hecke algebras over $\mathbb{Z}[q,q^{-1}]$ with $q$ an
indeterminant. However, it is not easy to write the basis
explicitly. Then one should ask, is there any basis which can be
written explicitly? In \cite{F}, Francis gave an integral
¡°minimal¡± basis for the center of a Hecke algebra. Then in
\cite{F2}, he used the minimal basis approach to provide a way of
describing and calculating elements of the minimal basis for the
center of an Iwahori-Hecke algebra which is entirely combinatorial.
In \cite{FJ2}, Francis and Jones found an explicit non-recursive
expression for the coefficients appearing in these linear
combinations for the Hecke algebras of type A. The relations between
the so-called Jucys-Murphy elements and centers of Hecke algebras are
also investigated widely. In \cite{DG2}, Dipper and James conjectured
that the center of a Hecke algebra of type A consists of symmetric
polynomials in the Jucys-Murphy elements. This conjecture was proved by
Francis and Graham \cite{FG} in 2006. An analogous conjecture for
Ariki-Koike algebras is still open.

The fact that Hecke algebras of finite type are all cellular leads
us to considering how to describe the centers of Hecke algebras by
cellular bases. Furthermore, how to describe the center of a
cellular algebra in general? Clearly, most of the approaches for
studying Hecke algebras can not be used directly for cellular
algebras, since we have no Weyl group structure to use. Then we must
look for some new method. In fact, the symmetry of Hecke algebras
provides us a way. We will do some work on the centers of symmetric
cellular algebras in this paper.

In order to describe our result exactly, we fix some notations
first. Let $A$ be a symmetric cellular $R$-algebra with a
non-degenerate symmetric bilinear form $f: A\times A\rightarrow R$.
Then $f$ determines a map $\tau: A\rightarrow R$ which is defined by
$\tau(a)=f(a,1)$ for every $a\in A$. We call $\tau$ a
symmetrizing trace. Denote by $\{D_{S,T}^\lam \mid S,T\in
M(\lam),\lam\in\Lambda\}$ the dual basis determined by $\tau$.  Let
$H(A)=\{\sum\limits_{\lam\in\Lambda, S,T\in M(\lam)}C_{S,T}^\lam
aD_{S,T}^\lam\mid a\in A\}$. It is the Higman ideal of $Z(A)$. For
any $\lam\in\Lambda$ and some $T\in M(\lam)$, set
$x_{\lam}=\sum\limits_{S\in M(\lam)} C_{S,T}^\lam D_{S,T}^\lam$, where $x_{\lam}$ is independent of $T$, and
$L(A)=\{\sum\limits_{\lam\in\Lambda}r_{\lam}x_{\lam}\mid r_{\lam}\in
R\}$. Then the main result of this paper is as follows.\\

\noindent\emph{{\bf Theorem.} Let $A$ be a symmetric cellular
algebra with a cellular basis $\{C_{S,T}^\lam \mid S,T\in
M(\lam),\lam\in\Lambda\}$ and the dual basis $\{D_{S,T}^\lam \mid
S,T\in M(\lam),\lam\in\Lambda\}$ determined by a symmetrizing trace
$\tau$.
Then}\\
\noindent (1)\, \emph{$L(A)$ is an ideal of $Z(A)$ and contains the
Higman ideal.}

\noindent (2)\, \emph{$L(A)$ is independent of the choice of
$\tau$}.

\noindent (3)\, \emph{If $R$ is a field, then the dimension of
$L(A)$ is not less than the number of non-isomorphic simple
$A$-modules.}

This theorem enlarge the well known Higman ideal to a new one for
the center of a symmetric cellular algebra.

\bigskip

\section{Preliminaries}

In this section, we first recall some basic results on symmetric
algebras and cellular algebras, which is needed in the following sections.
The so-called Higman ideal is also described. References for this
section are the books \cite{CR} and \cite{GP}.

Let $R$ be a commutative ring with identity and $A$ an associative
$R$-algebra. As an $R$-module, $A$ is finitely generated and free. Suppose that there exists an $R$-bilinear map
$f:A\times A\rightarrow R$. We say that $f$ is
 non-degenerate if the determinant of the matrix
$(f(a_{i},a_{j}))_{a_{i},a_{j}\in B}$ is a unit in $R$ for some
$R$-basis $B$ of $A$. We say $f$ is associative if $f(ab,c)=f(a,bc)$
for all $a,b,c\in A$, and symmetric if $f(a,b)=f(b,a)$ for all
$a,b\in A$.
\begin{dfn}
An $R$-algebra $A$ is called symmetric if there is a non-degenerate
associative symmetric bilinear form $f$ on $A$. Define an $R$-linear map $\tau: A\rightarrow R$ by $\tau(a)=f(a,1)$.
We call $\tau$ a symmetrizing trace.
\end{dfn}

Let $A$ be a symmetric algebra with a basis $B=\{a_{i}\mid
i=1,\ldots,n\}$ and $\tau$ a symmetrizing trace. Denote by
$D=\{D_{i}\mid i=i,\ldots,n\}$ the basis determined by the
requirement that $\tau(D_{j}a_{i})=\delta_{ij}$ for all $i,
j=1,\ldots,n$. We will call $D$ the dual basis of $B$. For arbitrary $1\leq i,j \leq n$, we write
$a_{i}a_{j}=\sum\limits_{k}r_{ijk}a_{k}$, where $r_{ijk}\in R$. Fix
a $\tau$ for $A$. Then in \cite{L}, we proved the following lemma.
\begin{lem}\label{lmp1}
Let $A$ be a symmetric algebra with a basis $B$ and the dual
basis $D$. Then the following hold:
$a_{i}D_{j}=\sum\limits_{k}r_{kij}D_{k};\,\,\,\,\,D_{i}a_{j}=\sum\limits_{k}r_{jki}D_{k}.$\qed
\end{lem}

We now consider the set $\{\sum\limits_{i}D_{i}aa_{i}\mid a\in A\}$. It is well known
that $\{\sum\limits_{i}D_{i}aa_{i}\mid a\in A\}$ is an ideal of the
center of $A$, see \cite{CR}. We here give a direct proof by the
lemma above.
\begin{prop}
Let $A$ be a symmetric algebra with a basis $B$ and the dual
basis $D$. Then $\{\sum\limits_{i}D_{i}aa_{i}\mid a\in A\}$ is an
ideal of the center of $A$.
\end{prop}
\noindent {Proof:}\, For arbitrary $a_{j}\in B$ and $a\in A$, we
have
$$\sum_{i}D_{i}aa_{i}a_{j}=\sum_{i,k}r_{ijk}D_{i}aa_{k},$$ and
$$\sum_{i}a_{j}D_{i}aa_{i}=\sum_{i,k}r_{kji}D_{k}aa_{i}.$$
Obviously, the right sides of the above two equations are equal.
Then $\{\sum\limits_{i}D_{i}aa_{i}\mid a\in A\}\subseteq Z(A)$. It
is clear that the set is an ideal of the center of $A$.\qed

The ideal $H(A)=\{\sum\limits_{i}D_{i}aa_{i}\mid a\in A\}$ is called
Higman ideal of the center of the algebra $A$. It is independent the choice of the dual bases.

The following proposition is also proved in \cite{L}.
\begin{prop}
\label{propp}Suppose that $A$ is a symmetric $R$-algebra with a
basis $\{a_{i}\mid i=1, \cdots, n\}$.  Let $\tau, \tau^{'}$ be two
symmetrizing traces. Denote by $\{D_{i}\mid i=1, \cdots, n\}$ and $\{D_{i}^{'}\mid i=1,
\cdots, n\}$ the dual bases determined by $\tau$ and $\tau^{'}$ respectively. Then for
any $1\leq i \leq n$, we have
$D_{i}^{'}=\sum\limits_{j=1}^{n}\tau(a_{j}D_{i}^{'})D_{j}.$\qed
\end{prop}

We now recall the definition of cellular algebras
introduced by Graham and Lehrer \cite{GL} and some well known
results.
\begin{dfn} {\rm (\cite{GL} 1.1)}
\label{dfgl1}Let $R$ be a commutative ring with identity. An
associative unital $R$-algebra is called a cellular algebra with
cell datum $(\Lambda, M, C, i)$ if the following conditions are
satisfied:

(C1) The finite set $\Lambda$ is a poset. Associated with each
$\lam\in\Lambda$, there is a finite set $M(\lam)$. The algebra $A$
has an $R$-basis $\{C_{S,T}^\lam \mid S,T\in
M(\lam),\lam\in\Lambda\}$.

(C2) The map $i$ is an $R$-linear anti-automorphism of $A$ with
$i^{2}=id$ which sends $C_{S,T}^\lam$ to $C_{T,S}^\lam$.

(C3) If $\lam\in\Lambda$ and $S,T\in M(\lam)$, then for any element
$a\in A$, we have\\
$$aC_{S,T}^\lam\equiv\sum_{S^{'}\in
M(\lam)}r_{a}(S^{'},S)C_{S^{'},T}^{\lam} \,\,\,\,(\rm {mod}\,\,\,
 A(<\lam)),$$ where $r_{a}(S^{'},S)\in R$ is independent of $T$ and
where $A(<\lam)$ is the $R$-submodule of $A$ generated by
$\{C_{S^{''},T^{''}}^\mu \mid S^{''},T^{''}\in M(\mu),\mu<\lam\}$.

Apply $i$ to the equation in (C3), we obtain

$(C3^{'})\,\, C_{T,S}^\lam i(a)\equiv\sum\limits_{S^{'}\in
M(\lam)}r_{a}(S^{'},S)C_{T,S^{'}}^{\lam} \,\,\,\,(\rm mod
\,\,\,A(<\lam)).$
\end{dfn}

By Definition \ref{dfgl1}, it is easy to check that $$C_{S,S}^{\lam}C_{T,T}^{\lam}\equiv \Phi(S,T)C_{S,T}^{\lam}\quad(\rm mod
\,\,\,A(<\lam)),$$ where $\Phi(S,T)\in R$ depends only on $S$ and $T$.

Let $A$ be a cellular algebra with cell datum $(\Lambda, M, C, i)$.
We recall the definition of cell modules.
\begin{dfn} {\rm (\cite{GL} 2.1)}
\label{dfgl2}For each $\lam\in\Lambda$, define the left $A$-module
$W(\lam)$ as follows: $W(\lam)$ is a free $R$-module with basis
$\{C_{S}\mid S\in M(\lam)\}$ and $A$-action
defined by\\
$$aC_{S}=\sum_{S^{'}\in M(\lam)}r_{a}(S^{'},S)C_{S^{'}}
\,\,\,\,(a\in A,S\in M(\lam)),$$ where $r_{a}(S^{'},S)$ is the
element of $R$ defined in (C3).
\end{dfn}

For a cell module $W(\lam)$, define a bilinear form $\Phi
_{\lam}:\,\,W(\lam)\times W(\lam)\longrightarrow R$ by $\Phi
_{\lam}(C_{S},C_{T})=\Phi(S,T)$, extended bilinearly and define
$$\rad(\lam):= \{x\in W(\lam)\mid \Phi_{\lam}(x,y)=0
\,\,\,\text{for all} \,\,\,y\in W(\lam)\}.$$ Then Graham and Lehrer
proved the following results in \cite{GL}.

\begin{thm} {\rm\cite{GL}}
Let $K$ be a field and $A$ a finite dimensional cellular algebra.
Denote the $A$-module $W(\lam)/\rad \lam$ by $L_{\lam}$, where
$\lam\in\Lambda$ with $\Phi_{\lam}\neq 0$. Let
$\Lambda_{0}=\{\lam\in\Lambda\mid \Phi_{\lam}\neq 0\}$. Then the set
$\{L_{\lam}\mid \lam\in\Lambda_{0}\}$ is a complete set of
(representative of equivalence classes of ) absolutely simple
$A$-modules.\qed
\end{thm}

\bigskip

\section{Centers of symmetric cellular algebras}

Let $A$ be a symmetric cellular algebra with a cell datum $(\Lambda,
M, C, i)$. Denote the dual basis by $D=\{D_{S,T}^\lam \mid
S,T\in M(\lam),\lam\in\Lambda\}$, which satisfies
$$
\tau(C_{S,T}^{\lam}D_{U,V}^{\mu})=\begin{cases} 1,&\text{$\lam=\mu,\,\,\,S=U,\,\,\,T=V$;}\\
0,\, &\text{otherwise.}\end{cases}
$$

For any $\lam, \mu\in \Lambda$, $S,T\in M(\lam)$, $U,V\in M(\mu)$,
write
$$C_{S,T}^{\lam}C_{U,V}^{\mu}=\sum\limits_{\epsilon\in\Lambda,X,Y\in M(\epsilon)}
r_{(S,T,\lam),(U,V,\mu),(X,Y,\epsilon)}C_{X,Y}^{\epsilon}.$$

Then in \cite{L}, we proved the following lemma.
\begin{lem}\label{corc1}
Let $A$ be a symmetric cellular algebra with a basis $B$. Let $D$ be the dual basis determined by a
given $\tau$. For arbitrary $\lam,\mu\in\Lambda$ and $S,T,P,Q\in M(\lam)$, $U,V\in M(\mu)$, the following hold:\\
(1)\,\,\,\,$D_{U,V}^{\mu}C_{S,T}^{\lam}=\sum\limits_{\epsilon\in
\Lambda, X,Y\in
M(\epsilon)}r_{(S,T,\lam),(X,Y,\epsilon),(U,V,\mu)}D_{X,Y}^{\epsilon}.$\\
(2)\,\,\,\,$C_{S,T}^{\lam}D_{U,V}^{\mu}=\sum\limits_{\epsilon\in
\Lambda, X,Y\in
M(\epsilon)}r_{(X,Y,\epsilon),(S,T,\lam),(U,V,\mu)}D_{X,Y}^{\epsilon}.$\\
(3)\,\,\,\,$C_{S,T}^{\lam}D_{S,T}^{\lam}=C_{S,P}^{\lam}D_{S,P}^{\lam}.$\quad\quad
(4)\,\,\,\,$D_{S,T}^{\lam}C_{S,T}^{\lam}=D_{P,T}^{\lam}C_{P,T}^{\lam}.$\\
(5)\,\,\,\,$C_{S,T}^{\lam}D_{P,Q}^{\lam}=0\,\, if \,\,T\neq Q.$\quad
(6)\,\,\,\,$D_{P,Q}^{\lam}C_{S,T}^{\lam}=0\,\, if \,\,P\neq S.$\\
(7)\,\,\,\,$C_{S,T}^{\lam}D_{U,V}^{\mu}=0 \,\,\,\,if\,\,\, \mu\nleq \lam.$\quad
(8)\,\,\,\,$D_{U,V}^{\mu}C_{S,T}^{\lam}=0 \,\,\,\,if \,\,\,\mu\nleq
\lam.$\qed
\end{lem}

Let $A$ be a symmetric cellular $R$-algebra with a symmetrizing trace
$\tau$. The dual basis $\{D_{S,T}^\lam \mid S,T\in
M(\lam),\lam\in\Lambda\}$ is determined by $\tau$. Then the Higman ideal is
$H(A)=\{\sum\limits_{\lam\in\Lambda, S,T\in M(\lam)}C_{S,T}^\lam
aD_{S,T}^\lam\mid a\in A\}$. For
any $\lam\in\Lambda$ and $T\in M(\lam)$, set
$x_{\lam}=\sum\limits_{S\in M(\lam)} C_{S,T}^\lam D_{S,T}^\lam$ and
$L(A)=\{\sum\limits_{\lam\in\Lambda}r_{\lam}x_{\lam}\mid r_{\lam}\in
R\}$. Now we are in a position to give the main result of this paper.
\begin{thm} Let $A$ be a symmetric cellular
algebra with a cellular basis $\{C_{S,T}^\lam \mid S,T\in
M(\lam),\lam\in\Lambda\}$ and the dual basis $\{D_{S,T}^\lam \mid
S,T\in M(\lam),\lam\in\Lambda\}$ determined by a symmetrizing trace
$\tau$.
Then\\
\noindent (1)\, $L(A)$ is an ideal of $Z(A)$ and contains the
Higman ideal $H(A)$.

\noindent (2)\, $L(A)$ is independent of the choice of
$\tau$.

\noindent (3)\, If $R$ is a field, then the dimension of
$L(A)$ is not less than the number of non-isomorphic simple
$A$-modules.
\end{thm}
\noindent {Proof:}\, (1) Firstly, we show that $H(A)\subseteq L(A)$.

Clearly, we only need to show that
$l:=\sum\limits_{S,T\in
M(\lam),\lam\in\Lambda}C_{S,T}^{\lam}C_{U,V}^{\mu}D_{S,T}^{\lam}\in
L(A)$ for any $C_{U,V}^{\mu}\in B$, where $\mu\in\Lambda, U,V\in M(\mu)$. We divide $l$ into three parts:
$l=l_{\lam=\mu}+l_{\lam<\mu}+l_{\lam\nleq\mu}$, where
$$l_{\lam=\mu}:=\sum_{S,T\in M(\lam),\lam=\mu}C_{S,T}^{\lam}C_{U,V}^{\mu}D_{S,T}^{\lam}$$ and the other
two parts are defined similarly.

By Lemma \ref{corc1} (7),
$l_{\lam\nleq\mu}=0$.

We now show that
$$l_{\lam=\mu}=\Phi_{\mu}(C_{U},C_{V})x_{\mu}.$$

By Lemma \ref{corc1} (5), $C_{U,V}^{\mu}D_{X,Y}^{\mu}=0$ if $V\neq Y$.
Then $l_{\lam=\mu}=\sum\limits_{X\in
M(\mu)}C_{X,V}^{\mu}C_{U,V}^{\mu}D_{X,V}^{\mu}$. By Definition \ref{dfgl1}, we have
$$l_{\lam=\mu}=\sum_{X\in
M(\mu)}\Phi_{\mu}(C_{U},C_{V})C_{X,V}^{\mu}D_{X,V}^{\mu}+\sum_{\eta<\mu, P,Q\in M(\eta)}r_{P,Q,\eta}C_{P,Q}^{\eta}D_{X,V}^{\mu},$$
where $r_{P,Q,\eta}\in R$. Note that by Lemma \ref{corc1} (7), $\sum\limits_{\eta<\mu, P,Q\in M(\mu)}r_{P,Q,\eta}C_{P,Q}^{\eta}D_{X,V}^{\mu}=0$, then
$l_{\lam=\mu}=\Phi_{\mu}(C_{U},C_{V})x_{\mu}$. This implies that $l_{\lam=\mu}\in L(A)$.

Now let us consider $l_{\lam<\mu}$. For arbitrary $\lam <\mu$, we
show that
$$\sum_{S,T\in M(\lam)
}C_{S,T}^{\lam}C_{U,V}^{\mu}D_{S,T}^{\lam}=\sum_{T\in
M(\lam)}r_{(S,T,\lam),(U,V,\mu),(S,T,\lam)}x_{\lam}.$$ Note that  $$\sum_{S,T\in M(\lam)
}C_{S,T}^{\lam}C_{U,V}^{\mu}D_{S,T}^{\lam}=\sum_{S,T\in
M(\lam)}(\sum_{\epsilon\in \Lambda, X,Y\in
M(\epsilon)}r_{(S,T,\lam),(U,V,\mu),(X,Y,\epsilon)}C_{X,Y}^{\epsilon})D_{S,T}^{\lam}.$$
 By $(C3)^{'}$ of Definition \ref{dfgl1},
if $\epsilon\nleq\lam$, then $r_{(S,T,\lam),(U,V,\mu),(X,Y,\epsilon)}=0$.
By Lemma \ref{corc1} (7), if $\epsilon<\lam$, then $C_{X,Y}^{\epsilon}D_{S,T}^{\lam}=0$. Thus
$$\sum_{S,T\in M(\lam)
}C_{S,T}^{\lam}C_{U,V}^{\mu}D_{S,T}^{\lam}=\sum_{S,T\in
M(\lam)}\sum_{X,Y\in
M(\lam)}r_{(S,T,\lam),(U,V,\mu),(X,Y,\lam)}C_{X,Y}^{\lam}D_{S,T}^{\lam}.$$
 By $(C3)^{'}$ of Definition \ref{dfgl1},
if $X\neq S$, then $r_{(S,T,\lam),(U,V,\mu),(X,Y,\lam)}=0$. By Lemma \ref{corc1} (5), if $Y\neq T$, then
$C_{X,Y}^{\lam}D_{S,T}^{\lam}=0$. Hence,
$$\sum_{S,T\in M(\lam)
}C_{S,T}^{\lam}C_{U,V}^{\mu}D_{S,T}^{\lam}=\sum_{S,T\in
M(\lam)}r_{(S,T,\lam),(U,V,\mu),(S,T,\lam)}C_{S,T}^{\lam}D_{S,T}^{\lam}.$$
Note that for arbitrary $S,S^{'}\in M(\lam)$, by $(C3)^{'}$ of Definition \ref{dfgl1}.
$$r_{(S,T,\lam),(U,V,\mu),(S,T,\lam)}=r_{(S^{'},T,\lam),(U,V,\mu),(S^{'},T,\lam)}.$$
We get
$$\sum_{S,T\in M(\lam)
}C_{S,T}^{\lam}C_{U,V}^{\mu}D_{S,T}^{\lam}=\sum_{T\in
M(\lam)}r_{(S,T,\lam),(U,V,\mu),(S,T,\lam)}x_{\lam}.$$ This implies
$l_{\lam<\mu}\in L(A)$. Then we obtain $l\in L(A)$.

\smallskip

Secondly, we show that $L(A)\subseteq Z(A).$

We only need to show that
$x_{\lam}C_{U,V}^{\mu}=C_{U,V}^{\mu}x_{\lam}$ for
arbitrary $\lam\in\Lambda$ and $\mu\in\Lambda, U,V\in M(\mu)$.

On one hand, by Lemma \ref{corc1} (1),
\begin{eqnarray*}
x_{\lam}C_{U,V}^{\mu}&=&\sum_{S\in
M(\lam)}C_{S,T}^{\lam}D_{S,T}^{\lam}C_{U,V}^{\mu}\\&=&\sum_{S\in
M(\lam)}\sum_{\epsilon\in\Lambda,X,Y\in
M(\epsilon)}r_{(U,V,\mu),(X,Y,\epsilon),(S,T,\lam)}C_{S,T}^{\lam}D_{X,Y}^{\epsilon}.
\end{eqnarray*}
By a similar method as in the first part, we get
$$x_{\lam}C_{U,V}^{\mu}=\sum_{S,X\in
M(\lam)}r_{(U,V,\mu),(X,T,\lam),(S,T,\lam)}C_{S,T}^{\lam}D_{X,T}^{\lam}.$$
On the other hand,
\begin{eqnarray*}
C_{U,V}^{\mu}x_{\lam}&=&\sum_{S\in
M(\lam)}\sum_{\epsilon\in\Lambda,X,Y\in
M(\epsilon)}r_{(U,V,\mu),(S,T,\lam),(X,Y,\epsilon)}C_{X,Y}^{\epsilon}D_{S,T}^{\lam}\\
&=&\sum_{S,X\in
M(\lam)}r_{(U,V,\mu),(S,T,\lam),(X,T,\lam)}C_{X,T}^{\lam}D_{S,T}^{\lam}.
\end{eqnarray*}
So $x_{\lam}C_{U,V}^{\mu}=C_{U,V}^{\mu}x_{\lam}$ for
arbitrary $\lam, \mu\in\Lambda$, $U,V\in M(\mu)$, that is, $L(A)\subseteq Z(A)$.

\smallskip

Finally, we show that $L(A)$ is an ideal of $Z(A)$.

It suffices to show that for arbitrary $c\in Z(A)$ and
$\lam\in\Lambda$, the element $cx_{\lam}\in L(A)$, that is,
$\sum\limits_{S\in M(\lam)}C_{S,T}^{\lam}cD_{S,T}^{\lam}\in L(A)$.

Since $c$ is $R$-linear combination of elements of $B$, then we only
need to prove that for arbitrary $C_{U,V}^{\mu}\in B$, the element
$\sum\limits_{S\in
M(\lam)}C_{S,T}^{\lam}C_{U,V}^{\mu}D_{S,T}^{\lam}\in L(A)$. Clearly,
this element is equal to $$\sum_{S\in
M(\lam)}\sum_{\epsilon\in\Lambda,X,Y\in
M(\epsilon)}r_{(S,T,\lam),(U,V,\mu),(X,Y,\epsilon)}C_{X,Y}^{\epsilon}D_{S,T}^{\lam}.$$
We know that it is equal to
$r_{(S,T,\lam),(U,V,\mu),(X,Y,\epsilon)}x_{\lam}$ by a similar way
as in the first part. This implies that
$\sum\limits_{S\in M(\lam)}C_{S,T}^{\lam}cD_{S,T}^{\lam}\in L(A)$.\\

(2)\,\,$L(A)$ is independent of the choice of $\tau$.

Let $\tau$, $\tau^{'}$ be two non-equal symmetrizing traces and $D$,
$d$ the dual bases determined by $\tau$ and $\tau^{'}$ respectively.
For arbitrary $d_{S,T}^{\lam}\in d$, by Proposition \ref{propp}, we have
$$d_{S,T}^{\lam}=\sum_{\varepsilon\in\Lambda, X,Y\in
M(\varepsilon)}\tau(C_{X,Y}^{\varepsilon}d_{S,T}^{\lam})D_{X,Y}^{\varepsilon}.$$
Then by Lemma \ref{corc1},
\begin{eqnarray*}
\sum_{S\in M(\lam)}C_{S,T}^{\lam}d_{S,T}^{\lam}&=&\sum_{S\in
M(\lam)}\sum_{\varepsilon\in\Lambda, X,Y\in
M(\varepsilon)}\tau(C_{X,Y}^{\varepsilon}d_{S,T}^{\lam})C_{S,T}^{\lam}D_{X,Y}^{\varepsilon}\\
&=&\sum_{S\in M(\lam)}\sum_{X\in
M(\lam)}\tau(C_{X,T}^{\lam}d_{S,T}^{\lam})C_{S,T}^{\lam}D_{X,T}^{\lam}.
\end{eqnarray*}

By the definition of $\tau$, we have
$\tau(C_{X,T}^{\lam}d_{S,T}^{\lam})=\tau(d_{S,T}^{\lam}C_{X,T}^{\lam})$.
Then by Lemma \ref{corc1},
$\tau(d_{S,T}^{\lam}C_{X,T}^{\lam})=0$ if $X\neq S$, that is,
$$\sum_{S\in M(\lam)}C_{S,T}^{\lam}d_{S,T}^{\lam}=\sum_{S\in
M(\lam)}\tau(C_{S,T}^{\lam}d_{S,T}^{\lam})C_{S,T}^{\lam}D_{S,T}^{\lam}.$$
We now  need to show  $\tau(C_{S,T}^{\lam}d_{S,T}^{\lam})$ is
independent of $S$. It is clear by the equations
$d_{S,T}^{\lam}C_{S,T}^{\lam}=d_{S^{'},T}^{\lam}C_{S^{'},T}^{\lam}$
for arbitrary $S^{'}\in M(\lam)$.\\

(3) We only need to find $|\Lambda_{0}|$ $R$-linear
independent elements in $L(A)$, where $|\Lambda_{0}|$ is the
number of the elements in $\Lambda_{0}$. By the definition of
$\Lambda_{0}$, for each $\lam\in\Lambda_{0}$, there exist $S,T\in
M(\lam)$, such that $\Phi_{\lam}(C_{S},C_{T})\neq 0$. Write
$x_{\lam}=\sum\limits_{U\in M(\lam)}C_{U,T}^{\lam}D_{U,T}^{\lam}$.
By Lemma \ref{corc1}, we know that the coefficient of
$D_{S,T}^{\lam}$ in the expansion of $C_{S,T}^{\lam}D_{S,T}^{\lam}$
is
$r_{(S,T,\lam),(S,T,\lam),(S,T,\lam)}=\Phi_{\lam}(C_{S},C_{T})\neq
0$ and is $0$ in the expansion of $C_{U,T}^{\lam}D_{U,T}^{\lam}$ for
any $U\neq S$. That is, the coefficient of $D_{S,T}^{\lam}$ in the
expansion of $x_{\lam}$ is not zero. We also know that the
coefficient of $D_{S,T}^{\lam}$ in the expansion of $x_{\mu}$ is
zero for any $\mu\nleq \lam$. Now let
$\sum\limits_{\lam\in\Lambda_{0}}r_{\lam}x_{\lam}=0$ and $\mu$ a
minimal element in $\Lambda_{0}$. Then $r_{\mu}$ must be zero. By
induction, we know that $r_{\lam}=0$ for each $\lam\in\Lambda_{0}$.
This implies that $\{x_{\lam}\mid\lam\in\Lambda_{0}\}$ is $R$-linear
independent. That is, $\dim_{R} L(A)$ is not less than the number of
(representatives of equivalence classes of) simple $A$-modules. \qed

By a similar way, we obtain the following result.

\begin{thm}Suppose that $R$ is a
commutative ring with identity and $A$ a symmetric cellular algebra
over $R$ with a cellular basis $B=\{C_{S,T}^\lam \mid S,T\in
M(\lam),\lam\in\Lambda\}$ and the dual basis $D=\{D_{S,T}^\lam \mid
S,T\in M(\lam),\lam\in\Lambda\}$ . Denote the set of the $R$-linear
combination of the elements of the set
$\{x_{\lam}^{'}=\sum\limits_{T\in
 M(\lam)}D_{S,T}^{\lam}C_{S,T}^{\lam}\mid \lam\in\Lambda\}$ by
 $L(A)^{'}$. Then the following hold:\\
(1) $L(A)^{'}$ is an ideal of $Z(A)$ and contains the Higman
ideal.\\
(2) $L(A)^{'}$ is independent of the choice of $\tau$.\\
(3) If $R$ is a field, then the dimension of $L(A)^{'}$ is not less
than the number of (representatives of equivalence classes of)
simple $A$-modules.
\end{thm}

We now give some examples of $L(A)$.\\
{\bf Example}\,\,\,Let $K$ be a field and $Q$
the following quiver\\
$$\xymatrix@C=13mm{
  \bullet \ar@<2.5pt>[r]^{\alpha_1}  & \bullet \ar@<2.5pt>[r]^(0.4){\alpha_2}
  \ar@<2.5pt>@[r][l]^(1){1}^{\alpha_1'}^(0){2}
  &\bullet\ar@<2.5pt>@[r][l]^(0.20){3}^(0.6){\alpha_2'}\cdots
  \bullet\ar@<2.5pt>[r]^(0.6){\alpha_{n-1}} & \bullet\ar@<2.5pt>@[r][l]^(0.35){\alpha_{n-1}'}^(0.75){n-1}^(0){n}\\
}$$ with relation $\rho$ given as follows.

(1) all paths of length $\geq 3$;

(2) $\alpha_{i}^{'}\alpha_{i}-\alpha_{i+1}\alpha_{i+1}^{'}$,
$i=1,\cdots, n-2$;

(3) $\alpha_{i}\alpha_{i+1}$, $\alpha_{i+1}^{'}\alpha_{i}^{'}$,
$i=1,\cdots, n-2$.

Let $A=K(Q,\rho)$. Define $\tau$ by

(1) $\tau(e_{1})=\cdots=\tau(e_{n})=1$;

(2)
$\tau(\alpha_{i}\alpha_{i}^{'})=\tau(\alpha_{i}^{'}\alpha_{i})=1$,
$i=1,\cdots, n-1$;

(3)$\tau(\alpha_{i})=\tau(\alpha_{i}^{'})=0$.

Then $A$ is a symmetric cellular algebra with a cellular basis

\[ \begin{matrix}
\begin{matrix} e_{1}\end{matrix} ;&
\begin{matrix} \alpha_{1}\alpha_{1}^{'} & \alpha_{1}\\ \alpha_{1}^{'} &
e_{2}\end{matrix} ; &
\begin{matrix} \alpha_{2}\alpha_{2}^{'} & \alpha_{2}\\ \alpha_{2}^{'} &
e_{3}\end{matrix} ;& \cdots ;&
\begin{matrix} \alpha_{n-1}\alpha_{n-1}^{'} & \alpha_{n-1}\\ \alpha_{n-1}^{'} &
e_{n}\end{matrix} ;&
\begin{matrix} \alpha_{n-1}^{'}\alpha_{n-1}\end{matrix}.
\end{matrix} \]
The dual basis is
\[ \begin{matrix}
\begin{matrix} \alpha_{1}\alpha_{1}^{'}\end{matrix} ;&
\begin{matrix} e_{1} & \alpha_{1}^{'}\\ \alpha_{1} &
\alpha_{1}^{'}\alpha_{1}\end{matrix} ;&
\begin{matrix} e_{2} & \alpha_{2}^{'}\\ \alpha_{2} &
\alpha_{2}^{'}\alpha_{2}\end{matrix} ;& \cdots ;&
\begin{matrix} e_{n-1} & \alpha_{n-1}^{'}\\ \alpha_{n-1} &
\alpha_{n-1}^{'}\alpha_{n-1}\end{matrix} ;&
\begin{matrix} e_{n}\end{matrix}.
\end{matrix} \]

It is easy to know that $L(A)$ is an ideal of $Z(A)$ generated by
$\{\alpha_{1}\alpha_{1}^{'},
\alpha_{1}\alpha_{1}^{'}+\alpha_{2}\alpha_{2}^{'},
\alpha_{2}\alpha_{2}^{'}+\alpha_{3}\alpha_{3}^{'}, \cdots,
\alpha_{n-2}\alpha_{n-2}^{'}+\alpha_{n-1}\alpha_{n-1}^{'},
\alpha_{n-1}^{'}\alpha_{n-1}\}$ and $H(A)$ is generated by
$\{2\alpha_{1}\alpha_{1}^{'}+\alpha_{2}\alpha_{2}^{'},
\alpha_{1}\alpha_{1}^{'}+2\alpha_{2}\alpha_{2}^{'}+\alpha_{3}\alpha_{3}^{'},
\alpha_{2}\alpha_{2}^{'}+2\alpha_{3}\alpha_{3}^{'}+\alpha_{4}\alpha_{4}^{'},
\cdots,
\alpha_{n-3}\alpha_{n-3}^{'}+2\alpha_{n-2}\alpha_{n-2}^{'}+\alpha_{n-1}\alpha_{n-1}^{'},
\alpha_{n-2}\alpha_{n-2}^{'}+2\alpha_{n-1}\alpha_{n-1}^{'}\}$.

Then $\dim_{K}L(A)=n$ since the rank of the matrix below is $n$.
$$\begin{bmatrix}
1 & 0 & 0 & 0 & \cdots & 0 & 0\\
1 & 1 & 0 & 0 & \cdots & 0 & 0\\
0 & 1 & 1 & 0 & \cdots & 0 & 0\\
& \cdots & & \cdots & & \cdots & & &\\
0 & 0 & 0 & 0 & \cdots & 1 & 1\\
0 & 0 & 0 & 0 & \cdots & 0 & 1
\end{bmatrix}_{(n+1)\times n}$$

We know that $\dim_{K}H(A)< n$ if $CharK$ is a factor of $n+1$ and
$\dim_{K}H(A)= n$ otherwise, since the determinant of the matrix
below is $n+1$.
$$\begin{bmatrix}
2 & 1 & 0 & 0 & \cdots & 0 & 0\\
1 & 2 & 1 & 0 & \cdots & 0 & 0\\
0 & 1 & 2 & 1 & \cdots & 0 & 0\\
& \cdots & & \cdots & & \cdots & & &\\
0 & 0 & 0 & 0 & \cdots & 2 & 1\\
0 & 0 & 0 & 0 & \cdots & 1 & 2
\end{bmatrix}_{n\times n}$$

Then $H(A)\subsetneq L(A)$ if $CharK$ is a factor of $n+1$ and
$H(A)= L(A)$ otherwise.

\smallskip

\noindent{\bf Example}\,\,\,Let $K$ be a field and $A=K[x]/(x^{n})$,
where $n\in \mathbb{N}$. Clearly, $A$ is a symmetric cellular
algebra with a basis $1, \bar{x}, \ldots, \overline{x^{n-1}}$. It is
easy to know that $L(A)$ has a basis $\overline{x^{n-1}}$ and
$Z(A)=A$. This example shows that $\dim_{K} Z(A)-\dim_{K} L(A)$ may
be very large.

\bigskip
\begin{prop}
Notations are as in Theorem A. Then  $x_{\lam}x_{\mu}=0$ for
arbitrary $\lam,\mu\in \Lambda$ with $\lam\neq\mu$.
\end{prop}
\noindent{Proof:}\, For arbitrary $\lam,\mu\in\Lambda$ with
$x_{\lam}x_{\mu}\neq 0$, then there exist $S_{0}\in M(\lam)$ and
$U_{0}\in M(\mu)$ such that
$$C_{S_{0},T}^{\lam}D_{S_{0},T}^{\lam}C_{U_{0},V}^{\mu}D_{U_{0},V}^{\mu}\neq
0.$$ This implies $D_{S_{0},T}^{\lam}C_{U_{0},V}^{\mu}\neq 0$. Then
by Corollary \ref{corc1}, there exists some $C_{X,Y}^{\epsilon}$
such that
$$r_{(U_{0},V,\mu),(X,Y,\epsilon),(S_{0},T,\lam)}\neq 0.$$ By (C3) of Definition \ref{dfgl1}, this implies that
$\lam\leq\mu$. By $x_{\lam}x_{\mu}=x_{\mu}x_{\lam}$, we get
$x_{\lam}x_{\mu}\neq 0$ also implies that $\mu\leq\lam$. Then $\lam
=\mu$ if $x_{\lam}x_{\mu}\neq 0$.\qed

\bigskip

\section{Semisimple case}

In this section, we consider the semisimple case. We will construct
all the central idempotents which are primitive in $Z(A)$ by
elements $x_{\lam}$ defined in Section 3 for a semisimple symmetric
cellular algebra.

Firstly, let us recall the definition of Schur elements. For
details, see \cite{GP}.

Let $R$ be a commutative ring with identity and $A$ an $R$-algebra.
Let $V$ be an $A$-module which is finitely generated and free over
$R$. The algebra homomorphism
$$\rho_{V}:A\rightarrow
 \rm End_{R}(V),\,\,\, \rho_{V}(a)v=av,\,\,\,\, {\text where}\,\,\,\,\, v\in V,\,\,\,\, a\in A,$$ is
called the representation afforded by $V$. The corresponding
character is the $R$-linear map defined by $$\chi_{V}:A\rightarrow
R,\,\,\,a\mapsto {\bf tr}(\rho_{V}(a)),$$ where {\bf tr} is the
usual trace of a matrix.

Let $K$ be a field and $A$ a finite dimensional symmetric
$K$-algebra with symmetrizing trace $\tau$. Let $B=\{B_{i}\mid i=1, \cdots, n\}$ be a basis and
$D=\{D_{i}\mid i=1, \cdots, n\}$ the dual basis determined by $\tau$. If $V$ is a split
simple $A$-module, denote the character by $\chi_{V}$, we have
$$\sum_{i}\chi_{V}(b_{i})\chi_{V}(D_{i})=c_{V}\dim_{K}V,$$ where $c_{V}\in K$ is the so-called Schur
element associated with $V$. We also denote it by $c_{\chi_V}$. It
is non-zero if and only if $V$ is a split simple projective
$A$-module \cite{GP}.
\begin{lem}
\label{lamgp}{\rm (\cite{GP} 7.2.7)} Let $A$ be a split semisimple
$K$-algebra. Then
$$\{e_{V}:=c_{V}^{-1}\sum_{i}\chi_{V}(b_{i})D_{i}\mid V\,\,\text{is a simple A-module}\}$$
is a complete set of central idempotents which are primitive in
$Z(A)$.
\end{lem}

Let $R$ be an integral domain and $A$ a symmetric cellular algebra
with cellular basis $\{C_{S,T}^\lam \mid S,T\in
M(\lam),\lam\in\Lambda\}$. Given a symmetrizing trace $\tau$, the
dual basis is $\{D_{S,T}^\lam \mid S,T\in M(\lam),\lam\in\Lambda\}$.
Let $K$ be the field of fractions of $R$. Define
$A_{K}:=A\bigotimes_{R}K$. Consider $A$ as a subalgebra of $A_{K}$
and extend $\tau$ of $A$ to $A_{K}$. Then we can construct all the
central idempotents which are primitive in $Z(A_{K})$ by $x_{\lam}$.
\begin{prop}
\label{propc3}  If $A_{K}$ is split semisimple, then
$\{c_{W(\lam)}^{-1}x_{\lam}\mid\lam\in\Lambda\}$ is a complete set
of central idempotents which are primitive in $Z(A_{K})$.
\end{prop}
\noindent {Proof:}\,\,The left $A_{K}$-module $W(\lam)$ is split
simple since $A_{K}$ is split semisimple. Then by Lemma \ref{lamgp},
we have
$$e_{W(\lam)}=c_{W(\lam)}^{-1}\sum_{\mu\in\Lambda, U,V\in M(\mu)}\chi_{W(\lam)}(C_{U,V}^{\mu})D_{U,V}^{\mu}.$$
Note that the character afforded by $W(\lam)$ is given by the
following
formula \\
$$\chi_{W(\lam)}(a)=\sum_{S\in M(\lam)}r_{a}(S,S)$$ for all $a\in
A$. Then we get $\chi_{W(\lam)}(C_{U,V}^{\mu})=\sum\limits_{S\in
M(\lam)}r_{(U,V,\mu),(S,T,\lam),(S,T,\lam)}.$ Then
\begin{eqnarray*}
e_{W(\lam)}&=&c_{W(\lam)}^{-1}\sum_{\mu\in\Lambda,U,V\in
M(\mu)}\sum_{S\in
M(\lam)}r_{(U,V,\mu),(S,T,\lam),(S,T,\lam)}D_{U,V,}^{\mu}\\
&=&c_{W(\lam)}^{-1}\sum_{S\in M(\lam)}\sum_{\mu\in\Lambda,U,V\in
M(\mu)}r_{(U,V,\mu),(S,T,\lam),(S,T,\lam)}D_{U,V,}^{\mu}\\
&=&c_{W(\lam)}^{-1}\sum_{S\in M(\lam)}C_{S,T}^{\lam}D_{S,T}^{\lam}.
\end{eqnarray*}\qed

\noindent{\bf Remark.}\, Clearly, $\{x_{\lam}\mid\lam\in\Lambda\}$
is a basis of the center of $A_{K}$ by this proposition. It is
different from the one in \cite{JY} when we consider Hecke algebras.

Note that the center of $A$ is equal to the intersection of $A$ and
the center of $A_{K}$. We now give a necessary condition for an
element of the center of $A_{K}$ being in $A$.
\begin{cor} Let $a_{\lam}\in K$ for all $\lam\in \Lambda$ and
 $a=\sum\limits_{\lam\in \Lambda}a_{\lam}x_{\lam}\in A$. Then
$a_{\lam}c_{W(\lam)}n_{\lam}\in R$ for arbitrary $\lam\in\Lambda$,
where $n_{\lam}$ is the number of elements in the set $M(\lam)$.
\end{cor}
\noindent {Proof:}\,\,For any $\lam\in\Lambda$, we know
$c_{W(\lam)}^{-1}x_{\lam}$ is a central idempotent of $A_{K}$ by
Proposition \ref{propc3}, i.e. $x_{\lam}^{2}=c_{W(\lam)}x_{\lam}$.
This implies that $ax_{\lam}=a_{\lam}c_{W(\lam)}x_{\lam}$. Clearly,
$ax_{\lam}\in A$ implies $\tau (ax_{\lam})\in R$. By the definition
of the dual basis, $\tau (x_{\lam})=m_{\lam}$. This completes the
proof.\qed

\bigskip\bigskip

\noindent{\bf Acknowledgments}\,\,\,The author acknowledges his supervisor Prof. C.C. Xi and the
support from the Research Fund of Doctor Program of Higher
Education, Ministry of Education of China. He also acknowledges Dr. Wei Hu for many helpful conversations.

\end{document}